\documentclass[12pt]{article}

\usepackage{amsthm}
\usepackage{amstext}
\usepackage{amssymb}
\usepackage{amsmath}
\usepackage{pb-diagram} \usepackage{amsthm}
\usepackage{rotating}
\usepackage{amsfonts}
\usepackage{graphicx}
\usepackage[
top    = 2cm,
bottom = 2cm,
left   = 2cm,
right  = 2cm]{geometry}
\theoremstyle{plain}
 \newtheorem{thm}{Theorem}
 \newtheorem{lem}{Lemma}
 
 \newtheorem{prop}{Proposition}
 
\theoremstyle{definition}

  \newtheorem*{thm*}{Theorem*}
 \newtheorem*{lem*}{Lemma*}
 \newtheorem*{cor*}{Corollary*}
 \newtheorem*{prop*}{Proposition*}
  \newtheorem*{defn*}{Definition}

\newcommand{\noi}{\noindent}

\begin{document}

\title{$2$-stratifolds  with fundamental group $\mathbb{Z}$}

\author{J. C. G\'{o}mez-Larra\~{n}aga\thanks{Centro de
Investigaci\'{o}n en Matem\'{a}ticas, A.P. 402, Guanajuato 36000, Gto. M\'{e}xico. jcarlos@cimat.mx} \and F.
Gonz\'alez-Acu\~na\thanks{Instituto de Matem\'aticas, UNAM, 62210 Cuernavaca, Morelos,
M\'{e}xico and Centro de
Investigaci\'{o}n en Matem\'{a}ticas, A.P. 402, Guanajuato 36000,
Gto. M\'{e}xico. fico@matem.unam.mx} \and Wolfgang
Heil\thanks{Department of Mathematics, Florida State University,
Tallahasee, FL 32306, USA. heil@math.fsu.edu}}
\date{}

\maketitle

\begin{abstract} $2$-stratifolds are a generalization of $2$-manifolds that occur as objects in applications such as in TDA. These spaces can be described by an associated bicoloured  labelled graph. In previous papers we obtained a classification of 1-connected trivalent $2$-stratifolds. In this paper we classify trivalent $2$-stratifolds that have fundamental group $\mathbb{Z}$. This classification implicitly  gives an efficient algorithm that applies to a bicoloured labelled graph to decide wether the associated $2$-stratifold has fundamental group $\mathbb{Z}$.
 \footnote {\bf{{AMS classification numbers}}:  57N10, 57M20, 57M05} 
\end{abstract}

%%%%%%%%%%%%%%%%%%%%%%%%%%%%%%%%%%%%%%%%
\section{Introduction}

A (closed) $2$-stratifold is a compact connected $2$-dimensional cell complex $X$ that can be constructed from a disjoint union of compact (connected) $2$-manifolds $W^2$ and a disjoint union $X^{(1)}$ of circles by attaching each component of $\partial W^2$ to $X^{(1)}$ via a covering map $\psi: \partial W^2 \to X^{(1)}$, with $\psi^{-1} (x)>2$ for $x\in X^{(1)}$. These $2$-stratifolds form a special class of $2$-dimensional stratified spaces. Matsuzaki and Ozawa \cite{MO} defined and studied a slightly more general class of $2$-dimensional stratified spaces, which they call {\it multibranched surfaces}, by allowing boundary curves, i.e. considering covering maps $\psi:\partial W'\to X^{(1)}$, where $\partial W'$ is a sub collection of the components of $\partial W^2$. A $2$-stratifold $X$ can be described by its bicoloured labelled graph $\Gamma_X$, whose white vertices are the components $W^2$, the black vertices are the components of $X^{(1)}$, and an edge is a component of $\partial W^2 \cap X^{(1)}$, where the label on an edge is the degree of the attaching map $\psi:\partial W^2 \to X^{(1)}$.\\

$2$-stratifolds arise as the nerve of certain Lusternick-Schnirelman type decompositions of $3$-manifolds into pieces where they determine whether the $\mathcal{G}$-category of the $3$-manifold is $2$ or $3$ (\cite{GGH}). They are related to {\it foams}, which include special spines of $3$-dimensional manifolds and which have been studied by Khovanov \cite{Ko} and Carter \cite{SC}. Matsuzaki and Ozawa showed that every $2$-stratifold $X$ embeds in $\mathbb{R}^4$, however $X$ embeds into some orientable closed $3$-manifold if and only if the singular set $X^{(1)}$ satisfies a certain regularity condition. In \cite{GGH3} we showed that very few $2$-stratifolds occur as spines of closed $3$-manifolds.\\

Graphs can be thought of as being $1$-stratifolds and persistent graphs occur in Topological Data Analysis \cite{VK}. The use of $2$-dimensional stratified spaces in TDA has been proposed in the literature \cite{B}, with a concrete example given in \cite{BH}. It seems likely that  $2$-stratifolds will be used in a similar fashion as models for applications.
For example they occur in the study of the energy landscape of cyclo-octane \cite{MT}, and $2$-stratifolds with boundary and foams with and without boundary
occur in the study of boundary singularities produced by the motion of soap films \cite{GT}.\\

        Closed $2$-manifolds are classified in terms of their fundamental groups, but for a given $2$-stratifold there are infinitely many non-homeomorphic $2$-stratifolds with the same fundamental group. Since the group can be computed from the labelled graph one would like to obtain conditions on the type of  bicoloured labelled graphs that correspond to a given fundamental group. As a first step in this direction we obtained in \cite{GGH1} conditions on the labelled graph $\Gamma_X$ such that the associated $2$-stratifold $X$ is $1$-connected and in \cite{GGH2} we give a complete classification of $1$-connected trivalent $2$-stratifolds in terms of their graphs. The term {\it trivalent} means that each point of $X^{(1)}$ has a neighborhood where three sheets meet.  For these $2$-stratifolds we developed an efficient algorithm to decide whether they are simply connected and implemented the algorithm in \cite{GGHH}.\\

         The ultimate goal is to obtain a classification of all $2$-stratifolds in terms of their labelled graphs.  In particular, with a view towards applications,  one would like to have efficient algorithms for deciding whether a given $2$-stratifold is of a certain type. In this paper we consider $2$-stratifolds $X$  whose fundamental group is infinite cyclic. 
In sections 3 and 4 we obtain necessary and sufficient conditions on $X$ and the type of the graph $\Gamma_X$ such that $\pi_1 (X)=\mathbb{Z}$.
However for a general $2$-stratifold graph $\Gamma_X$ it is difficult  to obtain conditions on the labelling that guarantee that $X$ has fundamental group $\mathbb{Z}$. We are able to do this in section 6 for trivalent $2$-stratifolds, after obtaining a classification for an interesting type of trivalent graphs, the echinus graphs, in section 5.   The main Theorem in section 6 provides a complete classification, in terms of conditions that can be read off from the labelled graph, of trivalent $2$-stratifolds with fundamental group $\mathbb{Z}$.\\

%%%%%%%%%%%%%%%%%%%%%%%%%%
\section{$2$-stratifolds and their graphs}

We first review the basic definitions and some results given in \cite{GGH1} and \cite{GGH2}. A  $2$-{\it stratifold} is a compact, Hausdorff space $X$ that contains a closed (possibly disconnected) $1$-manifold $X^{(1)}$ as a closed subspace with the following property: Each  point $x\in X^{(1)}$  has a neighborhood homeomorphic to $\mathbb{R}{\times}CF$, where $CF$ is the open cone on $F$ for some (finite) set $F$ of cardinality $>2$  and $X - X^{(1)}$ is a (possibly disconnected) $2$-manifold.\\

A component $B\approx S^1$ of $X^{(1)}$ has a regular neighborhood $N(B)= N_{\pi}(B)$ that is homeomorphic to $(Y {\times}[0,1]) /(y,1)\sim (h(y),0)$, where $Y$ is the closed cone on the discrete space $\{1,2,...,d\}$ and $h:Y\to Y$ is a homeomorphism whose restriction to $\{1,2,...,d\}$ is the permutation $\pi:\{1,2,...,d\}\to  \{1,2,...,d\}$. The space $N_{\pi}(B)$ depends only on the conjugacy class of $\pi \in S_d$ and therefore is determined by a partition of $d$. A component of $\partial N_{\pi}(B)$ corresponds to a summand of the partition determined by $\pi$. Here the neighborhoods $N(B)$ are chosen sufficiently small so that $N(B)$ is disjoint from $N(B' )$ for disjoint components $B$ and $B'$ of $X^{(1)}$. \\

For a given $2$- stratifold $X$ there is an associated bipartite graph $\Gamma=\Gamma_X$ embedded in $X$ as follows:\\

The white vertices $w_i$ of the graph $\Gamma$ are the components $W_i$ of the $2$-manifold $M:=\overline{X-\cup_j N(B_j)}$ where $B_j$  runs over the components of $X^{(1)}$. The black vertices $b_j$ are the $N(B_j)$'s. An edge $e_{ij}$ is a component $E_{ij}$ of $\partial M$; it joins a white vertex $w_i$ with a black vertex $b_j$  if $E_{ij}=W_i \cap N (B_j)$. Note that the number of boundary components of $W_i$ is the number of adjacent edges of $w_i$. 

We label the graph $\Gamma$ by assigning to a white vertex $w$ the genus $g$ of $W$ and by labelling an edge $e$ by $k$, where $k$ is the summand of the partition $\pi$ corresponding to the component $E$  of  $\partial N_{\pi}(B)$. (Here we use Neumann's \cite{N} convention of assigning negative genus $g$ to nonorientable surfaces; for example the genus $g$ of the projective plane or the Moebius band is $-1$, the genus of the Klein bottle is $-2$).  Note that the partition $\pi$ of a black vertex is determined by the labels of the adjacent edges.  \\

In sections 5 and 6 we obtain a classification of trivalent $2$-stratifolds in terms of their labelled graphs. 
Here a $2$-stratifold $X$ and its  labeled bicoloured graph $\Gamma_X$ are defined to be {\it trivalent}, if each black vertex $b$ is incident to either three edges each with label $1$ or to two edges, one with label $1$, the other with label $2$, or $b$ is a terminal vertex with adjacent edge of label $3$. This means that a neighborhood of a point of a component $B$ of the $1$-skeleton $X^{(1)}$ has $3$ sheets; the permutation $\pi:\{1,2,3\}\to  \{1,2,3\}$ of the regular neighborhood $N(B)= N_{\pi}(B)$ has partition $1+1+1$ or $1+2$ or $3$.\\

\noindent {\bf Notation}.  If $\Gamma$ is a bipartite labelled graph corresponding to the $2$-stratifold $X$ we let $X_{\Gamma} =X$ and $\Gamma_X =\Gamma$.  If $W$ is of genus $0$, we do not label $w$. \\

Another description of $X=X_\Gamma$ is as a quotient space $W \cup_{\psi}X^{(1)}$, where $W=\bigcup W_i$ and where $\psi:\partial W\to X^{(1)}$ is a covering map (and $|\psi^{-1}(x)| >2$ for every $x\in X^{(1)}$). For a component $E\approx S^1$ of $\partial W$ 
the label $m$ on the corresponding edge $e$ then corresponds to the attaching map $\psi ( z)=z^m$. An example is given in the picture below.\\
\begin{figure}[ht]
\begin{center}
\includegraphics[width=3.5in]{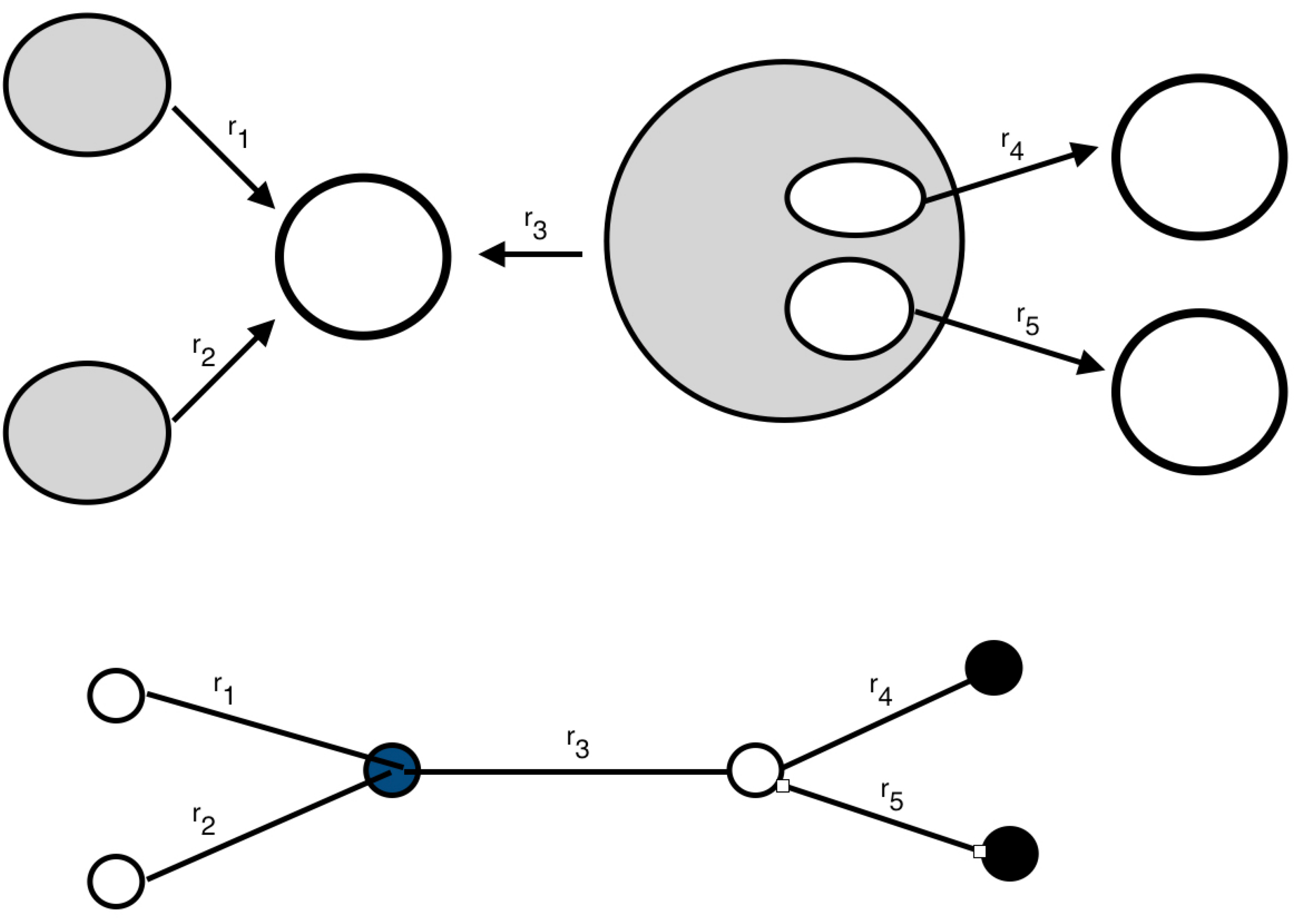} 
\end{center}

\end{figure}

If $\Gamma$ is a tree, then the labeled graph determines $X$ uniquely. If $\Gamma$ is not a tree then we need some additional data for $\Gamma$ which consists of a cocycle $\kappa \in H^1 (\Gamma; \mathbb{Z}_2) = Hom(H_1 (X), \mathbb{Z}_2 )$   together with a function $\lambda$ assigning to each edge $e$ of $\Gamma$ value $+1$ or $-1$  such that, for every simple cycle $c$ whose edges are $e_1,\dots,e_n$, $\lambda(e_1) \cdot .\,.\,.\, \cdot \lambda(e_n ) = \kappa([c])$. Here we identify $\mathbb{Z}_2 $, the cyclic group of order $2$, with $\{-1,1\}$ under multiplication.

One may assume such an evaluation $\lambda$ is constructed as follows: take a maximal tree  $T$ in $\Gamma$, and let $\lambda (e)=1$ if $e$ is an edge of $T$ and $\lambda (e)=\kappa([c])$ if $e$ is in $\Gamma - T$ and $c$ is the cycle of $T\cup e$. With this evaluation, the labelled graph $\Gamma$ together with a cocycle $\kappa$ determines $X_\Gamma$ uniquely.

 Thus if $\Gamma$ is homotopy equivalent to $S^1$  there is at most one (arbitrarily chosen) edge $e$ in the simple closed cycle $C$ of $\Gamma$ with $\lambda (e)=-1$.\\

The fundamental group of  $X_\Gamma$ can be computed from the labelled graph $\Gamma$. Generators of $\pi_1 (X_\Gamma )$ are the generators of $\pi_1 (W^2 )$ for each white vertex $w$, one generator for each black vertex $b$ (also denoted by $b$), and one generator $t$ for each edge outside a maximal tree $T$ of $\Gamma$. Relations are those of $\pi_1 (W^2 )$, $b^m =c$, for each edge $c \in T$ between $w$ and $b$ with label $m\geq 1$, and $t^{-1}c t=b^{\epsilon m}$, for each edge $c \in G-T$ between $w$ and $b$ with label $m$ ($\epsilon =\pm 1$). \\

In  \cite{GGH1} we showed the following.

\begin{prop}\label{retraction} There is a retraction $\rho:X_\Gamma  \to \Gamma_X$ such that $\rho^{-1}(b )$ is a singular curve $B\in X^{(1)}$ and $\rho^{-1}(w)$ is a $2$-manifold $W^2$.
\end{prop}

For simplicity, if $b$ is a black vertex  of $\Gamma_X$, we will sometimes write $b$ instead of $\rho^{-1} (b)$.  Also we will write $b$ instead of $[b]\in\pi_1 (X_\Gamma )$.\\

For a given  labelled graph $\Gamma$ we often prune away certain edges and vertices to obtain a subgraph $\Gamma_0$ such that there is an epimorphism $\pi_1 (X_\Gamma ) \to \pi_1 (X_{\Gamma_0})$ as follows: For a subgraph $\Gamma_0^{'}$ of $\Gamma$ let $Y'=\rho^{-1}(\Gamma_0^{'} )$. This is almost a $2$-stratifold, except that $Y'$ has possibly boundary curves corresponding to edges of $st(\Gamma_0^{'} )-\Gamma_0^{'}$, where $st(\Gamma_0^{'} )$ is the star of $\Gamma_0^{'}$ in $\Gamma$. Let $Y$ be the quotient of $X_\Gamma$ obtained by collapsing the closure of each component of $X_\Gamma  - Y'$ to a point i.e. $Y$ is obtained from $Y'$ by attaching disks to its boundary curves. Then $Y$ is a $2$-stratifold, $Y=X_{\Gamma_0}$, where $\Gamma_0$ is the union of $\Gamma_0^{'}$ and the labeled edges (with their vertices) of $st(\Gamma_0^{'} )-\Gamma_0^{'}$ which are adjacent to a black vertex of $\Gamma_0^{'}$. Then the quotient map $X_\Gamma \to X_{\Gamma_0}$ induces surjections of fundamental groups and first homology groups.\\

Using Proposition \ref{retraction} and this pruning construction we show in \cite{GGH1} the following

 \begin{prop}\label{simplyconnected} If $X$ is simply connected, then $\Gamma_X$ is a tree, all white vertices of $\Gamma_X$ have genus $0$, and all terminal vertices are white.
\end{prop}

%%%%%%%%%%%%%%%%%%%%%%%%%%%%%%%%%%%%%%%%
{\begin{section}{The homotopy type of the graph $\Gamma_X$}

In this section we show that a necessary condition for $\pi_1 (X_\Gamma )$ to be infinite cyclic is that the associated graph is homotopy equivalent to a circle,  the labels on all white vertices are $0$ and all terminal vertices are white.\\

 In the following Lemma $S_g$ denotes a closed surface of genus $g$.

\begin{lem}\label{circlegraph1}  Let $\Gamma=\Gamma_X$ be any graph of a $2$-stratifold $X=X_{\Gamma}$.\\
 (1) If $\Gamma$ has (at least) two black terminal vertices then there is an epimorphism $\pi_1 (X_{\Gamma})\twoheadrightarrow \mathbb{Z}_m *\mathbb{Z}_n$ for some $m,n\geq 3$.\\
 (2) If $\Gamma$ has a black terminal vertex and contains a white vertex of genus $g$ then there is an epimorphism $\pi_1 (X_{\Gamma})\twoheadrightarrow \mathbb{Z}_m *\pi_1 (S_g )$ for some $m\geq 3$.\\
(3) If $\Gamma$ contains two white vertices of genera $g_1$,  $g_2 $ then there is an epimorphism
$\pi_1 (X_{\Gamma})\twoheadrightarrow \pi_1 (S_{g{_1}}) *\pi_1 (S_{g_{2}} )$.
\end{lem}

\begin{proof} (1) Let $L$ be a linear subgraph of $\Gamma_X$ with terminal black vertices $b_1$, $b_2$ and prune $\Gamma$ at $L$. In the resulting graph construct $L_{mn}$ by deleting all edges (together with their white vertices) from each interior black vertex $b$, then adjoining to $b$ one edge with label $1$ (together with a white vertex) and finally changing the labels of all white vertices to $0$ (see Figure \ref{Zmn}). (This has the effect of killing $\rho^{-1}(b)$ in $\pi (X_\Gamma)$). Then there is a surjection $\pi_1 (X_{\Gamma})\twoheadrightarrow \pi_1 (X_{L_{mn}})\cong \mathbb{Z}_m * \mathbb{Z}_n$.\\

\includegraphics[width=6in]{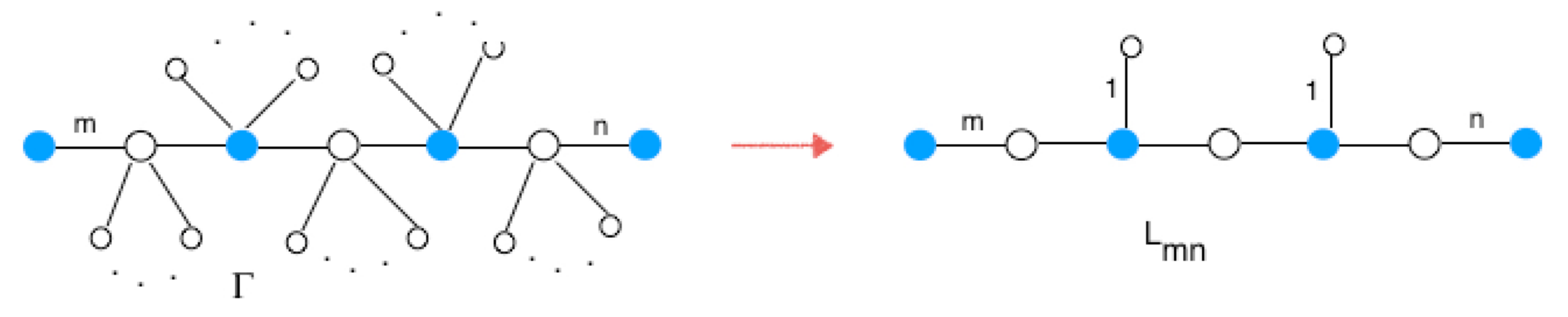}\label{Zmn} 

\noi(2) Let $L$ be a linear subgraph of $\Gamma_X$ with terminal black vertex $b$ and white vertex $w$ of genus $g$. Construct a graph $L_m$ as in (1) with terminal vertices $b$ and $w$. Then there is a surjection $\pi_1 (X_{\Gamma})\twoheadrightarrow \pi_1 (X_{L_m} )$. If $L_m$ is not  a one edge graph, then $\pi_1 (X_{L_{m}})\cong  \mathbb{Z}_m *\pi_1 (S_g )$, if $L_m$ is a one edge graph then  $\mathbb{Z}_m *\pi_1 (S_g )$ is a quotient of $\pi_1 (X_{L_m})$. 

\noi(3) In the proof of (2) replace  $L$ by a linear subgraph of $\Gamma_X$ with a terminal white vertex of genus $g_1$ and a terminal white vertex of genus $g_2$.
\end{proof}

\begin{lem}\label{finite} Let $\Gamma=\Gamma_X$ be a tree of a 2-stratifold $X_{\Gamma}$.\\
 (1) If $\Gamma_X$ has at most one black terminal vertex and all white vertices are of genus $0$ then $H_1 (X_{\Gamma})$ is  finite.\\
 (2) If $\Gamma_X$ has no black terminal vertices, contains at most one white vertex of genus $-1$ and all other white vertices are of genus $0$, then $H_1 (X_{\Gamma})$ is  finite.
\end{lem}

\begin{proof} (1)(a) Suppose $\Gamma$ has no black terminal vertices. The proof is by induction on the number of white vertices.
Let $w$ be a terminal (white) vertex with edge $e$ with label $m$ between $w$ and a black vertex $b$. If the degree of $b$ is $>2$ then $\Gamma_Y =\Gamma -(w\cup e)$ is as in the lemma and by induction $H_1 (Y)$ is finite. Therefore $H_1 (X)\cong H_1 (Y)/{\langle b^m \rangle}$ is finite.

\noi If the degree of $b$ is $2$ let $\Gamma_Y=\Gamma_X{-}(st(b)\cup\{w\}$. Then $H_1 (Y)\cong H_1 (X)/{\langle b \rangle}$, where
$b$ is a non-trivial element of the torsion subgroup (since $m\geq 1$). By induction $H_1 (Y)$ is finite and it follows that $H_1 (X)$ is finite.

\noi (1)(b) Now suppose $\Gamma_X$ has exactly one black terminal vertex $b$. If $b$ has finite order $n$ in $H_1 (X_\Gamma )$, let $\Gamma_Y$ be obtained from $\Gamma_X$ by attaching an edge $e$ with label $n$, together with a white vertex of genus $0$, to $b$. Then $\Gamma_Y$ is as in (1)(a) and so $H_1 (X_{\Gamma_X})=H_1 (X_{\Gamma_Y})$, is finite.

\noi We claim that $b$ does not have infinite order in $H_1 (X_\Gamma )$. To see this let $\Gamma_Y$ be the double of $\Gamma_X$ along $b$ (i.e. $\Gamma_Y$ is obtained from two disjoint copies of $\Gamma_X$ by identifying their copies of $b$). If $b$ has infinite order in $H_1 (X_\Gamma )$, then $\langle b\rangle$ is an infinite cyclic subgroup  of $H_1 (\Gamma_Y )=H_1 (X_\Gamma )\oplus H_1 (X_\Gamma )$ amalgamated along $\langle b\rangle$. This is a contradiction, since $\Gamma_Y$ is as in (1)(a).

(2) Let $\Gamma_{Y_{1}},\dots \Gamma_{Y_{r}}$ be the trees of $\Gamma_X{-}st(w)$ where $w$ is the white vertex of genus $-1$. By (1), $H_1 (\Gamma_{Y_{i}})$ is finite ($i=1,\dots ,r$).  Now $H_1 (X_{\Gamma})\cong H_1 (Y_1)\oplus \dots H_1 (Y_r ) \oplus \langle a \rangle$ /$\langle b_1^{m_1} \dots b_r ^{m_r} a^2 =1 \rangle$, where the $b_i$'s are the black vertices of $\overline{st(w)}$ and the $m_i$'s are the labels of the edges incident to $w$. It follows that $a^2 \in H_1 (Y_1)\oplus \dots H_1 (Y_r )$ and so $a$ has finite order in $H_1 (X_{\Gamma})$.

\end{proof}

\begin{lem}\label{circlegraph2} Let $\Gamma=\Gamma_X$ be a graph of a $2$-stratifold $X_{\Gamma}$ that is homotopy equivalent to $S^1$.\\
 (1) If $\Gamma_X$ has a white vertex of genus $g$ then there is an epimorphism $\pi_1 (X_{\Gamma})\twoheadrightarrow \mathbb{Z} * \pi_1 (S_g )$.\\
 (2) If $\Gamma_X$ has a black terminal vertex then there is an epimorphism $\pi_1 (X_{\Gamma})\twoheadrightarrow \mathbb{Z} * \mathbb{Z}_m$ for some $m\geq 3$. 
\end{lem}

\begin{proof} (1) If the white vertex $w$ of genus $g$ is a vertex of the cycle $C$ of $\Gamma$  let $C^*$ be the labelled graph obtained from $C$  by adjoining to each black vertex an edge (together with a white vertex) of label $1$ and assigning genus $0$ to all white vertices of $C^*$, except for $w$. Then there is an epimorphism $\pi_1 (X_{\Gamma})\twoheadrightarrow \pi (X_{C^*})= \mathbb{Z} * \pi_1 (S_g )$. 

If the white vertex $w$ of genus $g$ is not in $C$ connect $w$ to $C$ by a linear graph $L$ in $\Gamma$ such that $w$ is a terminal vertex of $L$ and $L\cap C$ is a vertex $v$. Let $C^*$ be as above an let $L^*$ be obtained from $C^* \cup L$ by adjoining to each black vertex of $L$ an edge (together with a white vertex) of label $1$ and assigning genus $0$ to all white vertices of $L^*$, except for $w$. Then there is an epimorphism $\pi_1 (X_{\Gamma})\twoheadrightarrow \pi (X_{L^*})= \mathbb{Z} * \pi_1 (S_g )$.

(2) Let $L$ be a linear subgraph of $\Gamma_X$ connecting the terminal black vertex $b$ (and with terminal edge labelled $m\geq 3$) to a  vertex of $C$ and construct $L^*$ as in (1). Then there is an epimorphism  $\pi_1 (X_{\Gamma})\twoheadrightarrow \pi (X_{L^*} )=\mathbb{Z} * \mathbb{Z}_m$. 
\end{proof}

\begin{prop}\label{circlegraph3} If $\pi_1 (X_{\Gamma})=\mathbb{Z}$ then  $\Gamma_X$ is homotopy equivalent to $S^1$, all white vertices of $\Gamma$ are of genus $0$, and all terminals vertices  are white.
\end{prop}

\begin{proof} The retraction $\rho :X_\Gamma\to\Gamma_X$ induces an epimorphism $\pi (X_\Gamma)\to\pi (\Gamma_X)$ an so $\Gamma_X$ is a tree or homotopy equivalent to $S^1$.   

\noi By Lemma \ref{circlegraph1}(1) $\Gamma_X$ has at most one black terminal vertex. \\
If $\Gamma_X$ has  one black terminal vertex then by Lemma \ref{circlegraph1}(2) all white vertices have genus $0$. Then by Lemma \ref{finite}(1) $\Gamma$ is not a tree.\\ 
If $\Gamma_X$ has  no black terminal vertices then by Lemma \ref{circlegraph1}(3) all white vertices have genus $0$, except possibly for one of genus $-1$. Then by Lemma \ref{finite}(2) $\Gamma$ is not a tree.\\
Hence $\Gamma$ is homotopy equivalent to $S^1$ and by Lemma \ref{circlegraph2} all white vertices have genus $0$ and all terminal vertices are white.
\end{proof}
\end{section}
%%%%%%%%%%%%%%%%%%%%%%%%%%%%%%%%%%%%%%%%%

\section{The general case of $\pi (X_\Gamma)\cong \mathbb{Z}$}

In this section we obtain the main result in the case when $\Gamma$ is any labelled graph satisfying the necessary conditions of Proposition \ref{circlegraph3}. By \cite{GGH2}, the requirements in Theorem \ref{main} below that certain stratifolds $Y_\Gamma$ be 1-connected can be determined from the labelled graphs $\Gamma_Y$.\\

We first obtain a Lemma that applies to any labelled graph $\Gamma$:

\begin{lem}\label{b=0} Let $\Gamma$ be a labelled graph and suppose $b$ is a black vertex of $\Gamma$ of degree $d\geq 2$ such that $\rho^{-1}(b)$ is contractible in $X_\Gamma$. Then $ \pi(X_\Gamma ) = \pi(X_{\Gamma_1})* ...* \pi(X_{\Gamma _c}) *F_r$ where $\Gamma_1,...,\Gamma_c$ are the components of $\Gamma_X{-}st\,b$ and $F_r$ is the free group of rank $r=d-c$.
\end{lem}

\begin{proof} Let $Y$ be the space obtained from $X_{\Gamma}$ by collapsing $\rho^{-1}(b )$ to a point $y_0$. Then there is an obvious map $p:Y\to \Gamma_X$ such that $p^{-1}(\Gamma_X -st\,b) = X_{\Gamma_1}\cup \dots \cup X_{\Gamma_c}$ and $p^{-1}(b)=y_0$. Since $\rho^{-1}(b )$ is contractible, $\pi (Y) \cong \pi (X_\Gamma )$. Now let $q: \Gamma_X \to \Gamma_Y$ be the quotient map, where $\Gamma_Y$ is obtained from $\Gamma_X$ by collapsing  each component $\Gamma_i$ to a point $y_i$ and let $T$ be a maximal tree in $\Gamma_Y$. Then $(qp)^{-1}(T)$ is homotopy equivalent to the wedge $X_{\Gamma_1}\vee\dots \vee X_{\Gamma _c}$. Since there are $r=d-c$ edges in $\Gamma_Y -T$ it follows that $Y=(qp)^{-1}(\Gamma_Y)$ is homotopy equivalent to $X_{\Gamma_1}\vee \dots \vee X_{\Gamma _c}\vee\bigvee_r S^1$.
\end{proof}

Recall from Proposition \ref{circlegraph3} that a necessary condition for  $\pi_1 (X_{\Gamma})=\mathbb{Z}$ is that $\Gamma_X$ is homotopy equivalent to $S^1$, all white vertices of $\Gamma$ are of genus $0$, and all terminals vertices  are white. The next lemma  lists more necessary conditions. By the {\it cycle} $C$ of $\Gamma$ we mean the subgraph that is homeomorphic to $S^1$.

\begin{lem}\label{necessary} Assume $\pi (X_\Gamma)\cong \mathbb{Z}$, then

\noi(a) $\Gamma$ is not homeomorphic to $S^1$.

\noi(b) There is at least one black vertex belonging to the cycle $C$ of $\Gamma$ whose degree is $>2$.

\noi(c) For each black vertex $b$, $\rho^{-1}(b)$ is homotopic to 0.

\noi(d) Let $B$ be a collection of black vertices containing at least one $b_1 \in C$. Then every component of $\Gamma{-}st B$ is 1-connected. 
\end{lem}

\begin{proof} We know that $\Gamma_X$ is homotopy equivalent to $S^1$, all white vertices of $\Gamma$ are of genus $0$, and all terminals vertices  are white.

\noi(a) We show that if $\Gamma_X =C$, a cycle,  then  $\pi(X_\Gamma)$ is infinite but not $\cong \mathbb{Z}$. 

Starting at a white vertex $w_1$ and reading the edges (counterclockwise) as\\
$w_1-{b}_1{-}w_2{-}b_2{-}\cdots{-}w_k{-} b_k{-}w_1$ with edge labels $m_1 , n_1 , \cdots , m_k ,n_k$ we obtain a presentation \\$\pi_1 (X_{\Gamma})=\{\,b_1 , \cdots , b_k , t\,\vert\, b_i^{ n_i}= b_{i+1}^{ m_{i+1}},  t^{-1}b_k^{n_k}t= b_1^{\epsilon m_1}\,,\,i=1,\dots,k{-}1\}$, where  $n_i$,$m_i \geq 1$ for $i=1,\dots,k$ and $\epsilon =\pm 1$. This group is an HNN extension of the nontrivial group (an iterated free product with amalgamation) 
$B := \{b_1,\dots, b_k \,\vert \, b_i^{n_i} = b_{i+1}^{m_{i+1}}\,,\,i=1,\dots,k-1\}$  and $b_1 ^{m_1}$ and $b_k^{n_k}$ are of infinite order in B.
Hence $\pi(X_\Gamma)\not\cong\mathbb{Z}$.

\noi(b) If each black vertex of $C$ has degree $2$ then pruning $\Gamma$ at each white vertex of $C$ yields a graph $\Gamma_C$ homeomorphic to $S^1$. Since $\rho_* :\pi (X_{\Gamma})\to \pi (X_{\Gamma_C})$ is surjective this contradicts  (a).

\noi(c) Since $\mathbb{Z}$ is Hopfian, the epimorphism $\rho_* : \pi(X_\Gamma) \to \pi(\Gamma)$ is an isomorphism and since $\rho( \rho^{-1}(b))$ is a point, the circle $ \rho^{-1}(b)$ is nulhomotopic in $X_\Gamma$.
%%%%%%

\noi(d) From (c) and Lemma \ref{b=0} it follows that $X_\Lambda$  is 1-connected for every component $\Lambda$ of $\Gamma{-}st(b_1 )$. If $\Lambda$ contains some $b\in B$, then since $\rho^{-1}(b)$ is contractible in $X_\Lambda$ it follows again from Lemma \ref{b=0} that every component of $\Gamma{-}st\{b_1,b\}$ is 1-connected. Now (d) follows from induction.\end{proof}

We remark that in (a) $H_1 (X_\Gamma) = \mathbb{Z}\oplus G_d$ where $G_d$ is a  group of order $d:=  n_1  n_2 \dots n_k - \epsilon m_1  m_2 ... m_k$,  the determinant of the $k$ by $k$ matrix representing $G_d =\{\,b_1 , \cdots , b_k , \,\vert\, [b_i ,b_j ]=1, b_i^{n_i}= b_{i+1}^{m_{i+1}} \,,b_k^{n_k}= b_1^{\epsilon m_1},\,i=1,\dots,k{-}1 \}$. If $X_\Gamma$ is trivalent (i.e. $0\leq m_i ,n_i \leq 2$, $m_i +n_i =3$) and $k > 1$ then $d$ is not $\pm 1$ and so $G_d \neq 0$. \\

The next proposition gives conditions that are necessary and sufficient.

\begin{prop}\label{mainprop} Let $\Gamma = \Gamma_X$ be a bicoloured labelled graph homotopy equivalent but not homeomorphic to $S^1$ and such that all white vertices have genus $0$ and all terminal vertices are white. Let $B$ be a collection of black vertices of $\Gamma$ containing a vertex $b_1$ of $C$ and such that $\rho^{-1}(b)$ is homotopic to $0$ in $X_\Gamma$ for each $b\in B$. Then $\pi (X_\Gamma)\cong \mathbb{Z}$ if and only all components of $X_{\Gamma -stB}$ are 1-connected.\end{prop}

\begin{proof} By Lemma \ref{necessary}(c)  the condition is necessary. Considering $b_1 \in B$ we have  by Lemma \ref{b=0} that $ \pi(X_\Gamma ) = \pi(X_{\Gamma_1})* ...* \pi(X_{\Gamma _c}) *\mathbb{Z}$, where $\Gamma_1,...,\Gamma_c$ are the components of $\Gamma_X{-}st\,b_1$. If $b\in B\cap \Gamma_i$ for some $1\leq i\leq n$, then again by Lemma \ref{b=0}, $\Gamma -stb$ splits $\pi (X_{\Gamma_i})$ into a free product. Doing this successively for all $b\in B$ we obtain $\pi(X_\Gamma )=(*\pi(X_{\Gamma_{ij}}))*\mathbb{Z}$, where the $\Gamma_{ij}$ are the components of $X-stB$. By assumption these are 1-connected and therefore $\pi(X_\Gamma )=\mathbb{Z}$ \end{proof}

We now state the main Theorem of this section.

\begin{thm}\label{main} Let $\Gamma$ be a bicoloured labelled graph. Then $\pi (X_\Gamma)\cong \mathbb{Z}$ if and only if the following four conditions hold.\\
(1) $\Gamma$ is homotopy equivalent but not homeomorphic to $S^1$.\\
(2) All white vertices have genus $0$ and all terminal vertices are white.\\
(3) The cycle of $\Gamma$ has a black vertex $b$ of degree $d > 2$ and $\rho^{-1}(b)$ is contractible in $X_\Gamma$.\\
(4) $X_{\Gamma_i}$  ($i = 1,\dots,d-1)$ is 1-connected where $\Gamma_1,\dots,\Gamma_{d-1}$ are the components of $\Gamma - st(b)$.
\end{thm}

\begin{proof} Proposition \ref{circlegraph3} and Lemma \ref{necessary}  show that the conditions are necessary.
Suppose (1)-(4) hold. Then for a black vertex $b$ as in (3) and (4) it follows from Lemma \ref{b=0} that $\pi(X_\Gamma) =\mathbb{ Z}$.
\end{proof}

\section{Trivalent echinus graphs}

In this section we consider bicoloured trivalent graphs $\Gamma =\Gamma_X$ that are homotopy equivalent to $S^1$ and such that the closure of each component of $\Gamma-C$ is a linear graph, where $C$ is the cycle of $\Gamma$.\\

First we list some necessary conditions for $\pi (X_\Gamma)\cong \mathbb{Z}$.

\begin{lem}\label{necessarytrivalent} Assume $\pi (X_\Gamma)\cong \mathbb{Z}$ and $\Gamma$ is trivalent, then

\noi(a) If all terminal edges of $\Gamma$ disjoint from $C$ have label $2$, then $\Gamma{-}C$ contains no white branch vertices.

\noi(b) Let $v$ be a non-terminal vertex of $\Gamma$ and let $L$ be linear graph that is the closure of a component of $\Gamma{-}\{v\}$ and contains a  terminal edge of $\Gamma$ of label $2$. If $v$ is a white vertex then $L$ is a linear $1212\dots 12$ graph. If $v$ is a black branch vertex then $L$ is a linear $11212\dots 12$ graph.

\noi(c) If all terminal edges of $\Gamma$ have label $2$ then $C$ contains no white branch vertices.
\end{lem}

\begin{proof} (a) Suppose $w$ is a white branch vertex in $\Gamma{-}C$ and $A_1 ,\dots, A_k$, $D$ are the closures of the components 
of $\Gamma{-}w$, where $D$ contains $C$. All terminal edges of $A_1\cup\dots\cup A_k$ have label $2$, and so 
$H_1 (X_{A_1\cup\dots\cup A_k})$ is not $0$. The retraction $\rho:X_D\to D$ shows that $H_1 (X_D )$ contains a $\mathbb{Z}$-summand and collapsing in $X_\Gamma$ the circle $\rho^{-1}w\cap \rho^{-1}(D)$ to a point we see that $H_1 (X_\Gamma)$ has $H_1 (X_D) + H_1 (X_{A_1\cup\dots\cup A_k})$ as a quotient, hence $H_1 (X_\Gamma)\not\cong \mathbb{Z}$.

\noi(b) Let $B$ be the set of black vertices of $\Gamma{-}L$. Then $L$ is a component of $\Gamma{-}B$ and is 1-connected by Lemma \ref{necessary}(d). 
Since all edge labels of $L$ are $1$ or $2$ (and in the second case the first label is $1$) the claim follows from Theorem 3 of [GGH].

\noi(c) By (a) we know that $\Gamma{-}C$ has no white branch vertices. If a component of $\Gamma{-}C$ contains a black branch vertex let $b$ be an outermost such vertex. Then the closures of the two components that do not meet $C$ are linear $1212\dots 12$ graphs of length $m_1 ,m_2$ with $m_2 \geq m_1$, say. Pruning away $L_2 \cup st(b)$ as in the figure below does not change the fundamental group. Hence we may assume that the closure of very component of $\Gamma{-}C$ is a linear graph as in (b).

 \includegraphics[width=4in]{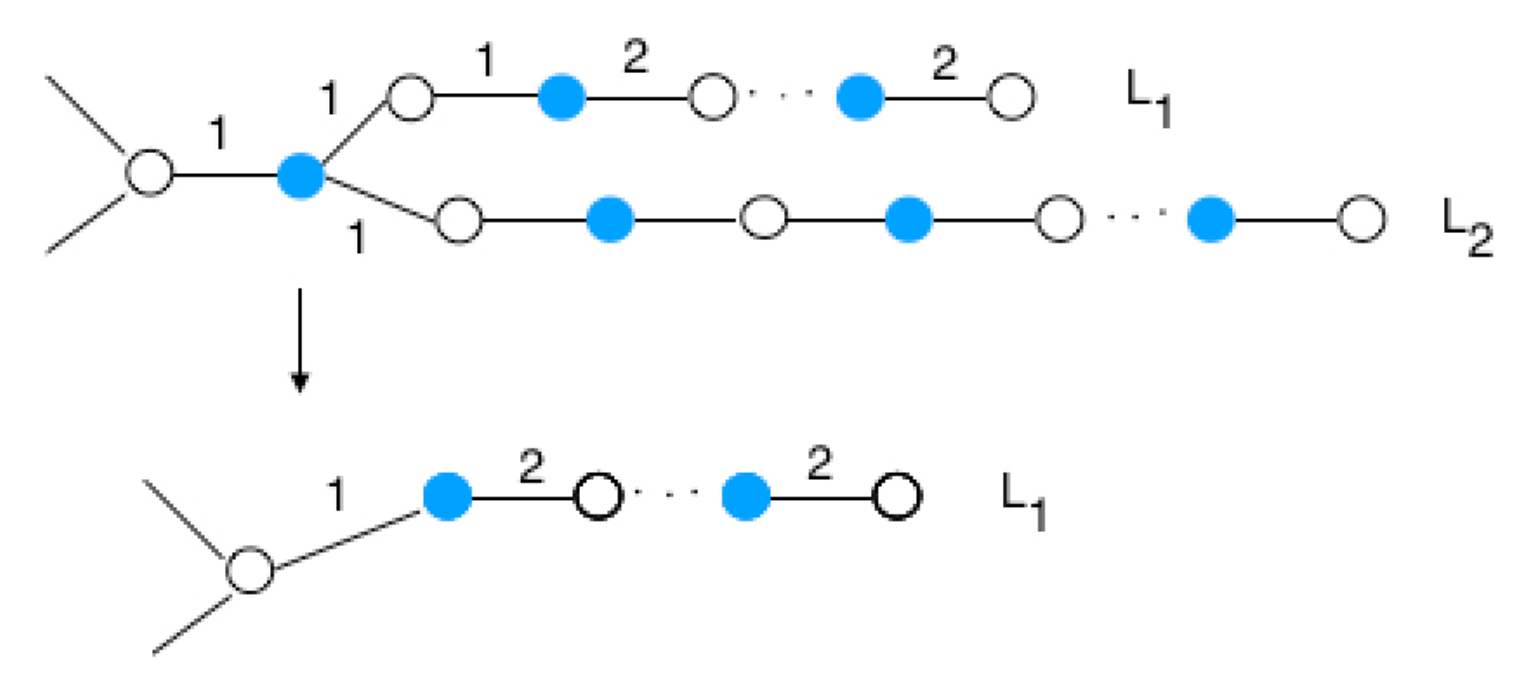} 

Let $w$ be a white vertex of $C$ of degree $>2$. We will show that $H_1 (X_\Gamma) \not\cong \mathbb{Z}$.\\
 If $w'$ is another such white vertex prune away all subgraphs other than $C$ that meet $w'$. For the resulting graph $\Gamma'$ there is an epimorphism $H_1 (X_\Gamma) \to H_1 (X_{\Gamma'})$. So may may assume that $C$ has exactly one vertex $w$ of degree $>2$.\\
 If degree$(w)>3$, prune away all subgraphs other than $C$ that meet $w$, except for one. Again there is an epimorphism $H_1 (X_\Gamma) \to H_1 (X_{\Gamma'})$. So may may assume that deg$(w)=3$.\\
 Note that the closure of each component of $\Gamma{-}C$ is a linear graph as in (2.2). If $L_r =b_1{-}w_{1}{-}b_{2}{-}w_{2}{-}\dots{-}b_r{-}w_r$ we may prune away the part after $w_{2}$ (the effect is killing $\rho^{-1}(b_1^2 )$ in $H_1$). So may may assume that all $r =2$.
 By the same argument, if $L=w{-}b_{1}{-}w_{1}{-}\dots{-}b_r{-}w_r$ we may (by pruning $L$ at $w_1$) assume that $L=w{-}b_{1}{-}w_{1}$.

If there is an edge $e$ of $C{-}st w$ with label $2$ let $w'$ be the white vertex incident with $e$ and let $\gamma$ be the path on $C$, containing $e$  with endpoints $w$ and $w'$. Then if $A$ is the union of $\gamma$ with the components of $\Gamma{-}C$ whose closure intersects $\gamma$, we have $H_1 (X_A )\neq 0$
since all the terminal edges of $A$ have label $2$. By killing all the loops in $X_\Gamma{-}X_A$ we see that $H_1 (X_\Gamma )$ has as quotient $\mathbb{Z} + H_1 (X_A )$ which is not $\mathbb{Z}$.

 If all edges of $C{-}stw$ have label 1 we compute $H_1 (X_\Gamma)$: Let $e_1$ and $e_2$ be the labels of the edges in $C$ adjacent to $w$ (the other edge $w{-}b_1$ adjacent to $w$ has label 1).  In $H_1 (X_\Gamma)$ one has $2 b_1 =0$ (where $b$ corresponds to a black vertex of $C$) and  $(e_2 -e_1)b = b_1$. Adjoining the relation $2b = 0$ (this is only needed if there is no black branch vertex) one obtains as a quotient of $H_1 (X_\Gamma)$ the abelian group
$( t, b, b_1 : 2b_1 =0, 2b=0, (e_2{-}e_1)b =b_1 )$  which is $\mathbb{Z} + \mathbb{Z}_2$ (no matter what $e_1$ and $e_2$ are).\end{proof}

A $\mathrm{p}$-{\it string} (of length $p$) is an oriented linear graph $w_0{-}b_1{-}w_1{-}b_2{-}\dots{-}b_k{-}w_k$ with all white vertices $w_i$ of genus $0$, successive edge labels $1212\dots 12$ (starting at $w_0$) and with $p$ labels of $2$.

A  $\mathrm{q}$-{\it string} is such an oriented linear graph with successive edge labels $2121\dots 21$ (starting at $w_0$) with $q$ labels of $2$.

If $L$ is a trivalent linear graph for which $X_L$ is 1-connected, then Theorem 3 of \cite{GGH} implies that $L$ does not contain a linear subgraph consisting of a $\mathrm{q}$-string followed by a $\mathrm{p}$-string for $q,p>0$. Thus we define $L(p,q)$ to be the linear graph consisting of a $\mathrm{p}$-string followed by a $\mathrm{q}$-string, i.e. $L(p,q)$ is the trivalent linear graph with terminal vertices white and edge labels $12\dots 1221\dots 21$, but we allow  $p$ or $q$ to be $0$.

We say that a $\mathrm{q}$-string is a {\it terminal} string of $\Gamma$  if $w_k$ is a terminal (white) vertex of $\Gamma$. Deleting from $\Gamma_X$ all terminal $\mathrm{q}$-strings, except for their initial vertices $w_0$, does not change the fundamental group of $X_\Gamma$. For the computation of the fundamental group it therefore  suffices to to restrict to {\it pruned} graphs, that are graphs without terminal $\mathrm{q}$-strings for $p>0$.\\

If $\Gamma$ is a trivalent pruned graph such that the closure of each component of $\Gamma-C$ is a linear graph and $\pi (X_\Gamma)=\mathbb{Z}$ then from Lemma \ref{necessarytrivalent}(c) we know that $C$ contains no white branch vertices and then it follows from Proposition \ref{main} and Theorem \ref{main} that all components of $\Gamma-stB$ are 1-connected $L(p,q)$-graphs, where $B$ is the set of black branch vertices of $C$. Furthermore, since $\Gamma$ is pruned, the closure of each component of $\Gamma-C$ is a linear graph $L(r,0)$-graph. We call these graphs {\it echinus}-graphs.

\begin{defn*} The echinus graph $E=E(p_1 ,q_1 ,r_1 ;\dots ;p_n ,q_n ,r_n )$ is the trivalent labelled graph $\Gamma$ with the following properties:\\
(1) $\Gamma$ is homotopy equivalent, but not homeomorphic to $S^1$.\\
(2) All vertices of $\Gamma$ are of degree $2$, except for $n\geq 1$ black vertices of the cycle $C$ of $\Gamma$.\\
(3) If $b_1 ,\dots,b_n$ are the successive black branch vertices for a fixed orientation of $C$ and $C_1 ,\dots,C_n $ are the successive components of $C{-}st\{b_1 ,\dots,b_n\}$, then $C_i$ is the linear graph $L(p_i ,q_i )$, with $p_i ,q_i \geq 0$, $i=1,\dots ,n$.\\
(4) The component of $\Gamma -st(b_i )$ that does not intersect $C$ is the linear graph $L(r_i ,0)$ with $r_i \geq 0$, $i=1,\dots ,n$.
\end{defn*}

 \includegraphics[width=6in]{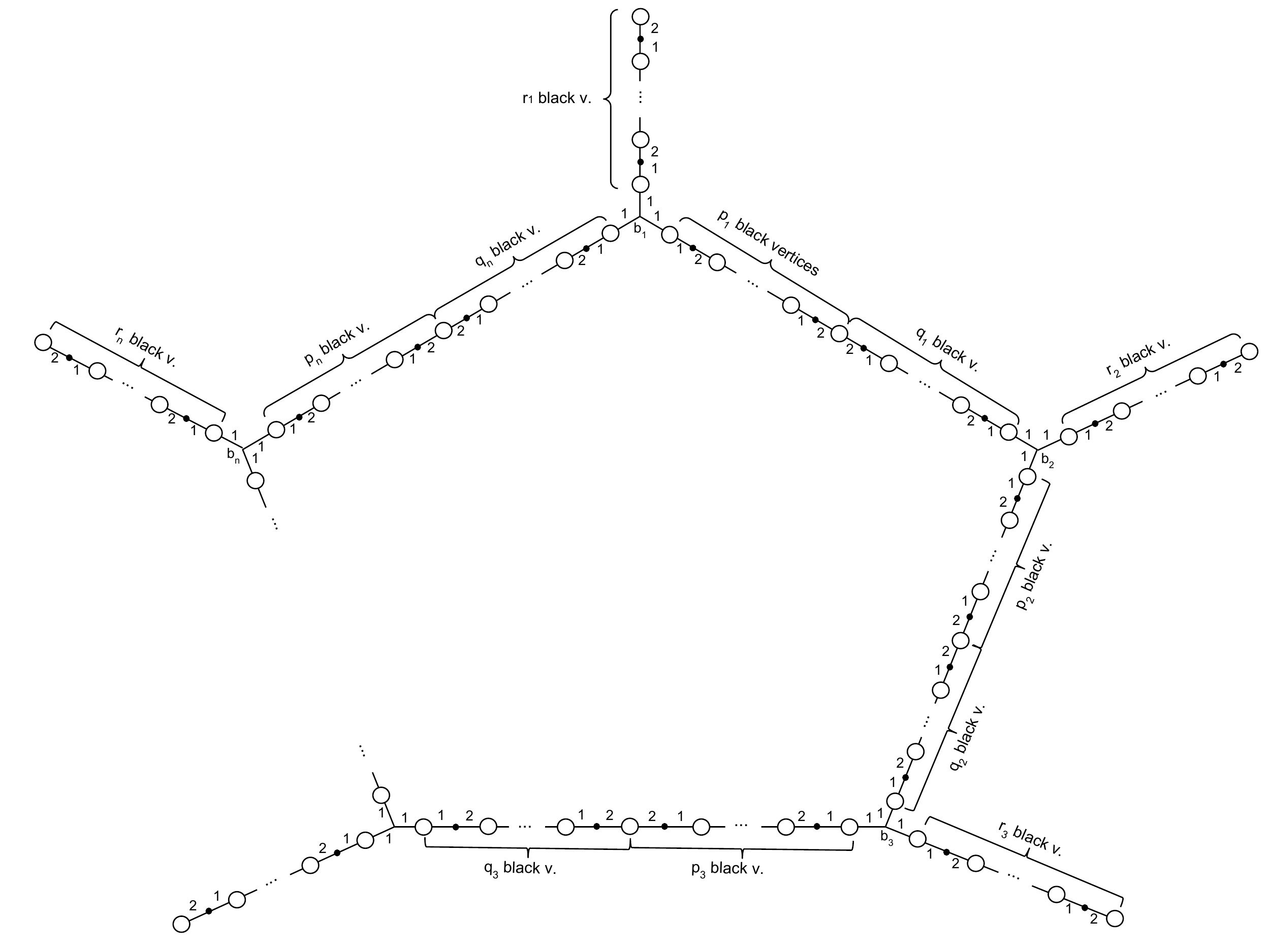}
 
The fundamental group of $X_E$ for $E=E(p_1 ,q_1 ,r_1 ;\dots ;p_n ,q_n ,r_n )$ has a presentation\\

$\pi (X_E )=\{b_1 ,\dots,b_n ,t\,:\,b_{j}^{2^{r_j}}=1, b_{i}^{2^{p_i}}=b_{i+1}^{2^{q_i}}, t^{-1} b_{n}^{2^{p_n}}t=b_{1}^{2^{q_1}},\,j=1,\dots,n,\,i=1,\dots,n{-}1\,\}$\\

and $H_1 (X_E )=\mathbb{Z}+T$, where $T$ has a relation matrix $M=\begin{pmatrix}
   A       \\
   B     
\end{pmatrix}$ as follows:\\
$A $ is the $n{\times}n$ diagonal matrix with the entry in the $i$th row and column  $2^{r_i}$;  $B$ is the $n{\times}n$ matrix ($b_{ij}$) where
$b_{ii} = 2^{p_i} , b_{ij} = -2^{q_i}$ if $j=i+1$ mod$n$ and all other $b_{ij}$'s are $0$. \\

We first consider echinus graphs that contain no branch vertices that are of distance $1$ to a terminal vertex (i.e. all $r_i >0$) and show:

\begin{prop}\label{H1E} Let $E=E(p_1 ,q_1 ,r_1 ;\dots ;p_n ,q_n ,r_n )$ with $r_i >0$ for all $i$.\\
Then $H_1 (X_E )=\mathbb{Z}$ if and only if exactly one of $\sum_{i=1}^n p_i$, $\sum_{i=1}^n q_i$ is $0$.
\end{prop}
\begin{proof} Recall that $T\cong\mathbb{Z}_{d_1}{\oplus}\dots\oplus\mathbb{Z}_{d_n}$ with $d_i \vert d_{i+1}$ and presentation matrix $D$ the diagonal matrix with $d_i$ the entry in the $i$th row and column. The $d_i$'s can be computed using elementary ideals as in \cite{CF} (Chapt.VIII Sect.4): $M$ is equivalent to $D$ and therefore
$(e_k) = E_k (M) = E_k (D) = (d_1 ,,d_2 ,\dots,d_{n-k})$, the subgroup generated by the product $d_1 d_2 \dots d_{n-k}$  if $0\leq k< n$, where $e_k$ ($\geq 0$) is the greatest common divisor of the determinants of all $(n-k){\times}(n-k)$ submatrices of M. Hence $d_i = e_{n - i} / e_{n - i +1}$ and in particular $T=0$ if and only if $e_0 = 1$.

For $M=\begin{pmatrix}A \\B\end{pmatrix}$ as above we have $detA=2^{r_1 + r_2 +\dots+ r_n}$ and one computes $detB = 2^{p_1 + p_2 +\dots+ p_n}-2^{q_1 + q_2 +\dots+q_n}$. All the other $n{\times}n$ submatrices contain a row of $A$ and therefore their determinants are (either $0$ or) divisible by $2$ (since by assumption each $r_i\neq 0$), and so $e_0 =1$ if and only if $2$ does not divide $2^{p_1 + p_2 +\dots+ p_n}-2^{q_1 + q_2 +\dots+q_n}$. This happens if and only if exactly one of $\sum_{i=1}^n p_i$, $\sum_{i=1}^n q_i$ is $0$.
\end{proof}

\begin{thm}\label{echinus1} Let $E = E(p_1, q_1, r_1; p_2, q_2, r_2; ... ; p_n, q_n, r_n)$ be an echinus graph with all $r_i \neq 0$.\\
Then the following are equivalent\\
(a) $\pi X_E = \mathbb{Z}$\\
(b) $H_1 (X_E ) = \mathbb{Z}$\\
(c) Precisely  one of $p_1 + p_2 + ... +p_n$ or $q_1 +q_2 +...+ q_n$  is 0. 
\end{thm} 

\begin{proof} By Proposition \ref{H1E} it suffices to show that if $q_1 +q_2 +...+ q_n =0$ but $p_1 + p_2 + ... +p_n\neq 0$ then $\pi X_E = \mathbb{Z}$.\\
 We may assume $r_1 = min\{r_1, ..., r_n\}$. Then eliminating $b_2,\dots,b_n$ in the presentation $of \pi X_E$ we obtain $\pi X_E = \{b, t : (b_1)^{2^{r_1}} = 1,  t^{-1}b_1^{2^m}t = b_1 \} $ where $m=p_1 + p_2 + ... +p_n >0$. \\
(*) If $r_1 \leq m$, then $b_1 = 1$ and so $\pi X_E = \{ t :  \} = \mathbb{Z}$.

If $r_1 > m$, raising the second relation to the power $2^{r_1  - m}$, we get \\
$\pi X_E = \pi X_E (p_1, q_1, r_1 - m ; p_2, q_2, r_2, ..., p_n, q_n, r_n)$.\\
After a number of these reductions of $r_1$ we obtain\\
$\pi X_E = \pi X_E (p_1, q_1, r'; p_2, q_2, r_2; ...; p_n, q_n, r_n )$ ,  where $r' \leq m$, which by (*) is $\mathbb{Z}$.\\
\end{proof}

We now consider echinus graphs that contain branch vertices of distance $1$ to terminal vertices. Note that if all $r_i$ are $0$, then $\pi (X_E )=1$.\\

Let $E$ be an echinus graph with cycle $C$ and let $B$ be the set of black branch vertices in $C$ at distance $1$ from a terminal vertex.
Assume $B$ is not empty and let $A$ be a component of $E{-}stB$. If $A$ does not consist of a single (white) vertex, $A$ is as in the figure below and with obvious notation we write \\

$A = A(p_1, q_1, r_1;  \dots; p_{n-1}, q_{n-1}, r_{n-1};p_n ,q_n)$, \qquad $n \geq 1, p_i \geq 0, q_i \geq 0, r_i > 0$.\\

 \includegraphics[width=5in]{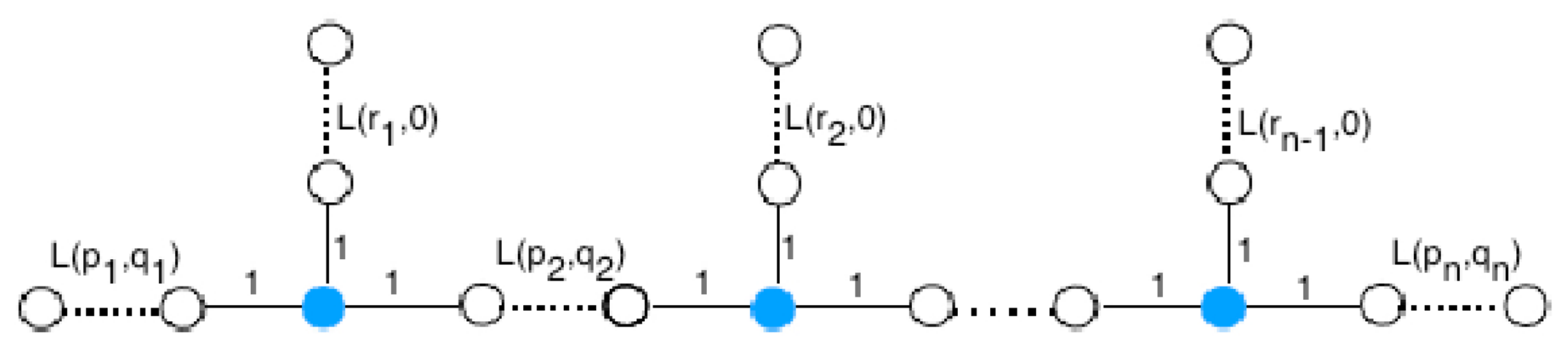}

\begin{prop}\label{echinusA} The following are equivalent :\\
(i) $X_A$ is 1-connected.\\
(ii) $H_1 (X_A)=0$.\\
(iii) If  $q_s >0$ for some $s<n$, then $p_t =0$ for all $t>s$.\end{prop}

\begin{proof} Clearly (i) implies (ii). We show that (ii) implies (iii). If (iii) does not hold, let $q_s >0$ and $p_t >0$ for some $s<t$ and let $w_1$ (resp. $w_2$) be the white starting vertex (resp. end vertex)
 of the $q_s$-string (resp. $p_t$-string). Consider the subgraph $A'$ of $A$ that is the union of the path $\gamma$ in $C\cap A$ from $w_1$ to $w_2$ and the components of $E{-}C$ whose closures intersect $\gamma$. Then  $\pi (X_{A'})$ is a quotient of $\pi (X_A)$  and $H_1 (X_{A'})\neq 0$ since all the terminal edges of $A'$ have label $2$ (recall all $r_i >0$). This contradicts (ii). Now we show that (iii) implies (1). This is clear if all $q_s =0$. If $q_s \neq 0$ then (iii) implies that $A = A(p_1, 0, r_1;  \dots; p_{s-1}, 0, r_{s-1}; p_s, q_s ,r_s ; 0, q_{s+1}, r_{s+1}; \dots ;  0,q_n)$. Pruning $A$ from the left we obtain a subgraph $A'=A(0, q_s ,r_s ; 0, q_{s+1}, r_{s+1}; \dots ;  0,q_n)$ with $\pi (X_A )=\pi (X_{A'})$ and then pruning $A'$ from the right we obtain $A''=w$, a single vertex, with $\pi (X_{A'} )=\pi (X_{A''})=1$.\end{proof}
 
 We now obtain a characterization of 1-connected echinus graphs that contain branch vertices that of distance $1$ to terminal vertices.
 
\begin{thm}\label{echinus2} Let $E$ be an echinus graph with cycle $C$ such that the set $B$  of black branch vertices in $C$ at distance $1$ from a terminal vertex is not empty. Then $X_E$ is 1-connected if and only if each component of $E{-}stB$ is a graph of type $A(p_1, 0, r_1;  \dots; p_{s-1}, 0, r_{s-1}; p_s, q_s ,r_s ; 0, q_{s+1}, r_{s+1}; \dots ;  0,q_n)$. \end{thm}

\begin{proof} This follows from Propositions \ref{mainprop} and \ref{echinusA}.
\end{proof}

%%%%%%%%%%
\section{Trivalent graphs}

In \cite{GGH2} we classified trivalent graphs that are 1-connected in terms of their graphs. We now describe necessary and sufficient conditions for  $\pi (X_\Gamma)$ to be $\mathbb{Z}$ that can be read off the graph for any trivalent graph $\Gamma$. 
We assume that $\Gamma$  satisfies the first three conditions of Theorem \ref{main} which can readily  be detected from the graph. Recall that {\it pruning} a graph means removing all terminal linear subgraphs with edge labels $21\dots 21$. This does not change the fundamental group of the associated $2$-stratifolds.  Thus we also assume that he graph $\Gamma$ is {\it pruned} i.e. there are no terminal linear subgraphs with edge labels $21\dots 21$. The To summarize, we assume\\
(1) $\Gamma$ is homotopy equivalent but not homeomorphic to $S^1$.\\
(2) All white vertices have genus $0$ and all terminal vertices are white.\\
(3) The cycle of $\Gamma$ has a black vertex $b$ of degree $d > 2$.\\
(4) $\Gamma$ is pruned.\\

For convenience we say that a labelled graph $\Delta$ is 1-connected if $X_{\Delta}$ is 1-connected. We now reduce $\Gamma$ to its core as follows.\\

If $\Gamma{-}C$ does not contain a black branch vertex of distance $1$ to a terminal vertex, we let $\Gamma_0 =\Gamma$ and say that $\Gamma_0$ is {\it core-reduced}. 
If $\Gamma{-}C$ contains  a black branch vertex of distance $1$ to a terminal vertex we let $B_0 =\{b_{01},\dots ,b_{0k}\}$ be the set of all such outermost black branch vertices, 
i.e. each $b_{0i}$ has distance $1$ from a terminal vertex $w_{0i}$ and the  component $T_{0i}$ of $\Gamma{-}(st(b_{0i} )\cup w_{0i})$ that does not meet $C$ does not contain a black branch vertex of distance $1$ to a terminal vertex. If some component $T_{0i}$  is not $1$-connected let $\Gamma_0 =\emptyset$.
If each $T_{0i}$  is 1-connected let $\Gamma'_0 =\Gamma{-}(st(B_0 )\cup \bigcup\{w_{0i}\} \cup \bigcup T_{0i} )$. If $\Gamma'_0 $ is pruned let $ \Gamma_0 =\Gamma'_0 $, otherwise let $\Gamma_0$ be the pruned $\Gamma'_0$.  Note that in this case $\pi (X_\Gamma )=\pi (X_{\Gamma_0})$, since $\rho^{-1}(b_{0i} )$ is homotopic to $0$ in $(X_\Gamma )$. 

Inductively, if $\Gamma_{n-1}{-}C$ contains an outermost black branch vertex of distance $1$ to a terminal vertex we let $B_{n-1} =\{b_{n-1,1},\dots ,b_{n-1,k_{n-1}}\}$ be the set of all such outermost black branch vertices that each has distance $1$ from a terminal vertex $w_{n-1,i}$  and let $T_{n-1,i}$ be the component of $\Gamma{-}(st(b_{n-1,i} )\cup w_{n-1,i})$. If some  $T_{n-1,i}$  is not $1$-connected we let $\Gamma_0 =\emptyset$.
If each $T_{n-1,i}$  is 1-connected we let $\Gamma_n$ be the pruned graph $\Gamma_{n-1}{-}(st(B_{n-1} )\cup \bigcup\{w_{n-1,i}\} \cup \bigcup T_{n-1,i} )$. 
We define the {\it core-reduced} subgraph $\Gamma_C$ as follows:

\noi $\Gamma_C =\begin{cases}
     \emptyset & \text{if } \Gamma_n =\emptyset \text{  for some } n\geq 0, \text{ otherwise}\\
     \Gamma_n & \text{ for the smallest $n$ such that } \Gamma_n{-}C \text{ does not contain a black branch vertex}\\
     & \text{ of distance $1$ to a terminal vertex}.
\end{cases}$

Note that if $\Gamma_C \neq \emptyset$ then $\pi (X_{\Gamma_C} )=\pi (X_{\Gamma})$.\\

We can now state our main classification Theorem.

\begin{thm}\label{classification} Let $\Gamma$ be a bicoloured pruned trivalent graph. Then  $\pi(X_\Gamma) = \mathbb{Z}$ if and only if\\
(a) The underlying graph of $\Gamma$ is homotopically equivalent, but not homeomorphic, to a circle\\
(b) All terminal vertices of $\Gamma$ are white and the genus of any white vertex is 0\\
(c) The core-subgraph $\Gamma_C$ is not the empty set\\
(d) If $B$ is the set of black vertices on the cycle $C$ of $\Gamma_C$ at distance 1 from a terminal vertex of $\Gamma_C$ then either\\
     (d1) $B$ is nonempty and $X_\Lambda$ is 1-connected for all components $\Lambda$ of  $\Gamma_C - st\,B$ or \\
     (d2) $B$ is empty, all branch vertices of $\Gamma_C$ are black, any edge of $\Gamma_C$ whose distance to $C$ (or, equivalently, to a terminal vertex of $\Gamma_C$) is odd has label 1, and $C$ is alternating. 
\end{thm}

Here the cycle $C$ of $\Gamma$ is alternating if, ignoring the edges incident to branch vertices, it is a $1212...12$ cycle with a positive (even) number of edges.

\begin{proof} 
Assume that $\pi(X_\Gamma) = \mathbb{Z}$. Then (a),(b),(c) and (d1)  follow from Theorem \ref{main} and Proposition \ref{mainprop}. As for (d2) we know from Lemma \ref{necessarytrivalent}(a) and (b)  that all branch vertices $\Gamma$ are black.  If there is an edge at odd distance $d$ from $C$ with label $2$ and $w$ is its white endpoint at distance $d$ from $C$ then, if $A$ is the closure of the  component of $\Gamma{-}w$, $X_A$ is not 1-connected because all terminal edges of $A$ have label $2$. This is a contradiction since $A$ occurs as a component as in Lemma \ref{necessary}(d).     

 By successively applying a pruning move as in the proof of of Lemma \ref{necessarytrivalent}(c)  (pruning off suitable terminal linear subgraphs that does not change the fundamental group), we may assume that the closures of all components of $\Gamma{-}C$ are disjoint linear subgraphs and, by further pruning, we may assume they have length $3$ with edge labels $1, 1, 2$. If the labelled graph is ``orientable" (i.e. all labels are positive) then it is an echinus graph $E = E(p_1, q_1, 1; p_2, q_2, 1;...; p_n,q_n,1)$ with $n>1$, $p_i \geq 0, q_i \geq 0$ and by Theorem \ref{echinus1} exactly one of $p_1+p_2+...+p_n$,   $q_1+q_2+...q_n = 0$. This means that $C$ is alternating.

If the labelled graph $E=E^-$  is not orientable then we may assume there is precisely one negative label in the graph and we may assume it is a $-1$ label of an edge of $C$ adjacent to (and ``to the left" of ) the black branch vertex $b_1$. Then $H_1 ( E^-)$  is presented by a $2n{\times}n$ matrix   $\begin{pmatrix} D\\F\\  \end{pmatrix}$
 where $D$ is an  $n{\times}n$ diagonal matrix with $2$'s on the diagonal and $F$ is an $n{\times}n$  matrix with determinant $2^{p_1+...+p_n} + 2^{q_1+... +q_n}$. Again, as in the proof of Proposition \ref{H1E}, $H_1 ( E^-)=\mathbb{Z}$, if and only if exactly one of $p_1+...+p_n$, $q_1+... +q_n$ is $0$, that is, $C$ is alternating.

To prove the converse, suppose (a), (b), (c) and (d) of the Theorem hold. If $b$ is an outermost black branch vertex of $\Gamma_C$ such that the closure of the two components 
of $\Gamma_C{-}b$ which do not intersect $C$ are linear $11212...12$ graphs $L_1, L_2$ of lengths $m_1, m_2$ respectively 
with $m_2 \geq m_1 \geq 1$, we prune off $L_2$ to obtain a labelled graph $\Gamma'_C$ with fewer such branch vertices. If $b$ is an outermost black branch vertex of $\Gamma_C$ adjacent to a terminal white vertex $w$ and $L\neq w$ is the closure of the component of $\Gamma_C{-}b$, delete $L\cup st(b)\cup w$ from $\Gamma_C$. If the resulting subgraph is not pruned, eliminate the terminal linear subgraph with edge labels $21\dots 21$ to get $\Gamma'$. In either case $\pi (X_{\Gamma_C})=\pi (X_{\Gamma'})$.
Repeating this pruning process we obtain an echinus graph $E$ with  $\pi (X_\Gamma)=\pi (X_{\Gamma_C})=\pi (X_E)$ and such that $E$ satisfies (a), (b), (c) and (d) of the Theorem. Now $\pi (X_E)=\mathbb{Z}$ by Theorems \ref{echinus1} and \ref{echinus2}.
\end{proof}

{\bf Acknowledgments:} J. C. G\'{o}mez-Larra\~{n}aga would like to thank INRIA Saclay, Francia and CNRS UMI 2001, Laboratorio Solomon Lefschetz, Mexico, for financial support.

\end{document}